\documentclass[10pt]{article}

\title{\Large \textbf{Mirror Descent and the Information Ratio}}
\author{Tor Lattimore and Andr\'as Gy\"orgy \\[0.1cm] {\small \tt \{lattimore,agyorgy\}@google.com}}
\date{DeepMind, London}

\usepackage{amsmath}
\usepackage{amssymb}
\usepackage{amsthm}
\usepackage{natbib}
\usepackage{bm}
\usepackage{mathrsfs}
\usepackage{mathtools}
\usepackage{nicefrac}
\usepackage[capitalise]{cleveref}
\usepackage[bf]{caption}
\usepackage{mdframed}
\usepackage{xspace}
\usepackage[boxed]{algorithm}
\usepackage[textsize=tiny]{todonotes}
\usepackage[cal=zapfc,bb=fourier]{mathalfa}

\theoremstyle{plain}
\newtheorem{theorem}{Theorem}

\newtheorem{lemma}[theorem]{Lemma}

\newtheorem{corollary}[theorem]{Corollary}

\theoremstyle{definition}
\newtheorem{definition}[theorem]{Definition}
\newtheorem{assumption}[theorem]{Assumption}

\newtheorem{remark}[theorem]{Remark}
\theoremstyle{remark}

\newcommand{\sP}{\mathcal P}
\newcommand{\sPp}{\mathcal P_{\scalebox{0.5}{$\bm{+}$}}}
\newcommand{\Reg}{\mathfrak R}
\newcommand{\BReg}{\mathfrak {BR}}

\newcommand{\ip}[1]{\langle #1 \rangle}

\newcommand{\cA}{\mathcal A}

\newcommand{\cV}{\mathcal V}
\newcommand{\cD}{\mathcal D}
\newcommand{\cX}{\mathcal X}
\newcommand{\cZ}{\mathcal Z}

\newcommand{\cI}{\mathfrak I}
\newcommand{\cG}{\mathcal G}

\newcommand{\epsd}{\epsilon_\cD}
\newcommand{\R}{\mathbb R}
\newcommand{\KL}{\operatorname{D}}
\newcommand{\norm}[1]{\Vert #1 \Vert}

\newcommand{\E}{\mathbb E}

\newcommand{\ones}{\bm{1}}
\newcommand{\bip}[1]{\left\langle #1 \right\rangle}
\newcommand{\argmin}{\operatornamewithlimits{arg\,min}}
\newcommand{\ceil}[1]{\lceil #1 \rceil}

\renewcommand{\d}[1]{\operatorname{d}\!#1}

\newcommand{\sind}{\bm{1}}
\newcommand{\diam}{\operatorname{diam}}
\newcommand{\dom}{\operatorname{dom}}
\newcommand{\zeros}{\bm{0}}
\newcommand{\prior}{A^*_{\text{\tiny{pr}}}}
\newcommand{\post}{A^*_{\text{\tiny{po}}}}

\newcommand{\relint}{\operatorname{relint}}
\newcommand{\conv}{\operatorname{conv}}

\begin{document}

\maketitle

\begin{abstract}
We establish a connection between the stability of mirror descent and the information ratio by \cite{RV14}.
Our analysis shows that mirror descent with suitable loss estimators and exploratory distributions enjoys the same bound on the adversarial regret as the bounds on the Bayesian regret
for information-directed sampling. 
Along the way, we develop the theory for information-directed sampling and provide an efficient algorithm for adversarial bandits for which the regret upper bound matches exactly the best known 
information-theoretic upper bound. 
\end{abstract}

\section{Introduction}

The combination of minimax duality and the information-theoretic machinery developed by \cite{RV14} has yielded a series of elegant arguments bounding the minimax regret
for a variety of regret minimisation problems. 
The downside is that the application of minimax duality makes the approach non-constructive. The \textit{existence} of certain policies is established without
identifying what those policies are. Our main contribution is to show that the information-theoretic machinery can be translated in a natural way to the language of online linear optimisation, yielding
explicit policies. Before you get too excited, these policies are not guaranteed to be efficient -- they must solve a convex optimisation problem that may be infinite dimensional. 
Nevertheless, it provides a clear path towards algorithm design and/or improved bounds, as we illustrate with an application to finite-armed bandits.

To maximise generality, our results are stated using the linear partial monitoring framework, which is flexible enough to model most classical setups.
Readers who are not familiar with partial monitoring should not be put off.
Our analysis does not depend on subtle concepts specific to finite partial monitoring, like the cell decomposition or observability. Examples
are given in \cref{tab:examples}.

A linear partial monitoring game is defined by an action space $\cA \subset \R^d$, a signal space $\Sigma$, a latent space $\cZ$ and two functions: a signal function $\Phi : \cA \times \cZ \to \Sigma$ 
and a loss function $\ell : \cZ \to [0,1]^d$. Both the signal and loss functions are known to the learner.
What is special about partial monitoring is that the learner never directly observes the realised losses, instead receiving signals that are correlated with the losses in a way that depends on the loss and signal functions. 
At the start of the game, an adversary secretly chooses a sequence $(z_t)_{t=1}^n$ with $z_t \in \cZ$. 
A policy is a mapping from action/signal sequences to distributions over actions.
The learner interacts with the environment over $n$ rounds.
In each round $t$, the learner uses their policy to find a distribution $P_t$ over the actions based on the history $(A_s)_{s=1}^{t-1}$ and $(\sigma_s)_{s=1}^{t-1}$, where $A_s$ is the action chosen in round $s$ and
$\sigma_s = \Phi_{A_s}(z_s)$ is the signal. The learner then samples $A_t$ from $P_t$ and observes the corresponding signal.
The regret of a policy $\pi$ is defined as
\begin{align*}
\Reg_n(\pi, (z_t)) = \max_{a \in \cA} \E\left[\sum_{t=1}^n \ip{A_t - a, \ell(z_t)}\right]\,,
\end{align*}
where the expectation integrates over the randomness in the actions chosen by the learner.
The arguments $\pi$ and $(z_t)$ are omitted when they are obvious from the context.
The quantity of interest is generally the minimax adversarial regret, defined as
\begin{align*}
\Reg_n^\star = \inf_{\pi} \sup_{(z_t)} \Reg_n(\pi, (z_t)) \,,
\end{align*}
where the infimum is taken over all policies of the learner and the supremum is over all possible choices of the adversary.
Given a finitely supported distribution $\mu$ on $\cZ^n$, the Bayesian regret of policy $\pi$ is
\begin{align*}
\BReg_n(\pi, \mu) = \int_{\cZ^n} \Reg_n(\pi, (z_t)) \d{\mu}((z_t))\,.
\end{align*}
A recently popular method for controlling the adversarial regret non-constructively appeals to minimax duality to show that
\begin{align}
\Reg_n^\star = \sup_{\mu} \inf_{\pi} \BReg_n(\pi, \mu)\,,
\label{eq:minimax}
\end{align}
where the supremum is over all finitely supported priors.
The Bayesian regret is then bounded uniformly over all priors using the information-theoretic argument of \cite{RV14}.
A limitation of this approach is that the application of minimax duality is non-constructive. It yields a bound
on the minimax regret but gives no hint towards an algorithm.

\paragraph{Contributions}
Our main contribution is a proof that bounds on the information ratio introduced by \cite{RV14} imply bounds on the stability of online learning algorithms mirror descent (MD) and
follow the regularised leader (FTRL).
The results provide an effortless proof of the main theorem of \cite{LS19pminfo} and hint towards the existence of improved algorithms for zeroth-order bandit convex optimisation.
Along the way, we further generalise the information-theoretic machinery to derive adaptive bounds and to make it more suitable for analysing games for which the minimax regret
is not $\Theta(n^{1/2})$.
A concrete consequence is an efficient algorithm for $d$-armed adversarial bandits for which $\Reg_n \leq \sqrt{2dn}$, improving on the best known result for an efficient algorithm that is
$\Reg_n \leq \sqrt{2dn} + 48k$ by \cite{ZL19}. A modest improvement that nevertheless illustrates the applicability of the approach.

\paragraph{Related work}
Mirror descent has its origins in the classical convex optimisation \citep{Nem79}, while follow the regularised leader goes back to the work by \cite{Gordon99}. 
As far as we know, the first application to bandits
was by \cite{AHR08}. The information-theoretic analysis for bandit problems was developed
in two influential papers by \cite{RV14,RV16}. These focussed on the Bayesian setting, with no connections made to the adversarial framework.
\cite{BDKP15} used minimax duality to argue that the minimax (adversarial) regret is equal to the worst-case Bayesian regret
and used this to derive the first proof that the minimax regret for convex bandits in one dimension is $O(\sqrt{n} \log(n))$.
The same plan has been used for convex bandits for larger dimensions \citep{BE18,Lat20-cvx} and partial monitoring \citep{LS19pminfo}, the latter of which establishes \cref{eq:minimax}
in the present setup. None of these works yields an efficient algorithm, but these have now been found for both settings \citep{BLE17,LS19pmsimple}, in both cases based on mirror descent.
Connections between the information ratio and mirror descent were investigated by \cite{ZL19}, who showed that bounds on the stability of mirror descent imply bounds on the
information ratio with somewhat restrictive assumptions. These results hinted at a deeper connection, but the analysis is somehow in the wrong direction, since the adversarial
regret is already a stronger notion than Bayesian regret.
The policy we propose in \cref{sec:exp-opt} is almost identical to the exploration by optimisation algorithm suggested by \cite{LS19pmsimple}. The difference is that now the bias of the loss
estimators is incorporated into the optimisation problem in a more natural way.

\section{Notation and conventions}\label{sec:notation}
Recall that a proper convex function $F : \R^d \to \R \cup \{\infty\}$ is Legendre if it is lower semi-continuous, essentially smooth and essentially strictly convex \citep[\S26]{Roc15}.
Throughout, let $F : \R^d \to \R \cup \{\infty\}$ be a Legendre function.

\begin{assumption}\label{eq:ass}
Let $\cD \subset \conv(\cA)$ be compact, convex and have non-empty relative interior, where $\conv(\cA)$ is the convex hull of $\cA$.
We make the following assumptions:
\begin{enumerate}
\item[\textit{(a)}] (finite action set): $1 < |\cA| < \infty$.
\item[\textit{(b)}] (bounded losses): $\ip{a, \ell(z)} \in [0,1]$ for all $a \in \cA$ and $z \in \cZ$.
\item[\textit{(c)}] (domain of potential): $\cD \subset \dom(F) \triangleq \{x \in \R^d : F(x) < \infty\}$.
\item[\textit{(d)}] (bounded potential): $\diam(\cD) = \sup_{x, y \in \cD} F(x) - F(y) < \infty$.
\end{enumerate}
\end{assumption}

The restriction to finite action sets avoids delicate measure-theoretic technicalities.
Note, since $\cD$ is compact, (d) is automatic when $F$ is continuous on $\cD$ with the subspace topology, which holds for all potentials considered
in the literature that satisfy (c).

\paragraph{Basic notation}
Precedence is given to the expectation operator: $\E[X]^\alpha$ denotes $(\E[X])^\alpha$ for random variables $X$ and reals $\alpha$.
The relative interior of a subset $A$ of a topological vector space is $\relint(A)$.
The standard basis vectors in $\R^d$ are $e_1,\ldots,e_d$.
Let $\sP$ be the space of probability distributions over $\cA$ and $\sPp = \{p \in \sP : p(a) > 0\,, \forall a \in \cA\}$ and $\sP_\epsilon = \{p \in \sP : p(a) \geq \epsilon\,, \forall a \in \cA\}$. Occasionally elements $p \in \sP$ are identified with vectors in $\R^d$ in the obvious way.
The Fenchel--Legendre dual of $F$ is the convex function defined by $F^\star(u) = \sup_{x \in \R^d} \ip{u, x} - F(x)$. Bregman divergences with respect to $F$ and $F^\star$ are
\begin{align*}
\KL(p, q) &= F(p) - F(q) - \ip{\nabla F(q), p - q} \\
\KL_\star(x, y) &= F^\star(x) - F^\star(y) - \ip{\nabla F^\star(y), x - y}\,.
\end{align*}
The assumption that $F$ is Legendre ensures that duality holds so that $(\nabla F)^{-1} = \nabla F^\star$ and
\begin{align}
\KL(p, q) = \KL_\star(\nabla F(q), \nabla F(p))\,.
\label{eq:kl-dual}
\end{align}
The space of finitely supported probability distributions on $\cZ \times \cD$ is denoted by $\cV$.
Finally, let 
\begin{align*}
\epsd = \max_{a \in \cA} \min_{b \in \cD} \max_{z \in \cZ} \ip{b - a, \ell(z)}\,,
\end{align*}
which vanishes in the typical case that $\cD = \conv(\cA)$.

\begin{table}[h!]
\centering
\renewcommand{\arraystretch}{1.6}
\scriptsize
\scalebox{0.95}{
\begin{tabular}{|llllll|}
\hline
\textsc{name} & $\cA$ & $\cZ$ & $\Sigma$ & $\ell(z)$ & $\Phi_a(z)$ \\ \hline
full information & $\{e_1,\ldots,e_d\}$ & $[0,1]^d$ & $[0,1]^d$ & $z$ & $z$ \\
$d$-armed bandits & $\{e_1,\ldots,e_d\}$ & $[0,1]^d$ & $[0,1]$ & $z$ & $z_a$ \\ 
linear bandits & arbitrary & $ \subset \R^d$ & $[0,1]$ & $z$ & $\ip{a, z}$ \\
graph feedback $(\dagger)$ & $\{e_1,\ldots,e_d\}$ & $[0,1]^d$ & $([d] \times [0,1])^*$ & $z$ & $(b, z_b)_{b \in N_a}$ \\
convex bandit $(\ddagger)$ & arbitrary & $\{\text{cvx } z \in [0,1]^\cA\}$ & $[0,1]$ & $(z(a))_{a \in \cA}$ & $z(a)$ \\
\hline
\multicolumn{6}{|p{11.8cm}|}{
$\dagger$ A bandit with graph feedback problem depends on a directed graph over the actions represented by a collection of sets $(N_a)_{a=1}^d$ with $N_a$
the set of edges originating from action $a$. When playing action $a$ the learner observes the losses for actions $b \in N_a$.
} \\ 
\multicolumn{6}{|p{11.8cm}|}{
$\ddagger$ The convex bandit problem is often formulated with a discrete action set. Alternatively, the first step in the analysis performs a discretisation.
} \\
\hline
\end{tabular}} 
\caption{Examples}\label{tab:examples}
\end{table}

\section{Mirror descent and FTRL}
Before presenting the new results, let us remind ourselves about the application of MD and FTRL to partial monitoring.
Given a sequence of loss estimates $(\hat \ell_t)_{t=1}^n$ with $\hat \ell_t \in \R^d$ and a sequence of non-increasing and strictly positive learning rates $(\eta_t)_{t=1}^n$,
MD produces a sequence $(q_t)_{t=1}^n$ with $q_t \in \cD$ defined inductively by
\begin{align}
q_1 &= \argmin_{q \in \cD} F(q) &
q_{t+1} &= \argmin_{q \in \cD} \ip{q, \hat \ell_t} + \frac{\KL(q, q_t)}{\eta_t} \,.
\label{eq:md}
\end{align}
Follow the regularised leader also produces a sequence $(q_t)_{t=1}^n$ with $q_t \in \cD$ defined by
\begin{align}
q_t = \argmin_{q \in \cD} \sum_{s=1}^{t-1} \ip{q, \hat \ell_s} + \frac{F(q)}{\eta_t}\,.
\label{eq:ftrl}
\end{align}
The next theorem bounds the regret of MD and FTRL with respect to the estimated losses.
There are many sources for results like this \citep[theorem 28.4, exercise 28.12]{LS20book}. 

\begin{theorem}\label{thm:md}
Suppose that one of the following is true: 
\begin{enumerate}
\item [(a)] $(q_t)_{t=1}^n$ are chosen according to \cref{eq:md} and $\eta_t = \eta$ is constant; or
\item [(b)] $(q_t)_{t=1}^n$ is chosen according to \cref{eq:ftrl}.
\end{enumerate}
Then,
\begin{align*}
\max_{a^* \in \cD} \sum_{t=1}^n \ip{q_t - a^*, \hat \ell_t} \leq \frac{\diam(\cD)}{\eta_n} + \sum_{t=1}^n \frac{\Psi_{q_t}(\eta_t \hat \ell_t)}{\eta_t} \,,
\end{align*}
where $\Psi_q(x) = \KL_\star(\nabla F(q) - x, \nabla F(q))$ is called the `variance' or `stability' term.
\end{theorem}

\begin{remark}
The function $x \mapsto \Psi_q(x)$ is convex and there is no randomness in any of the quantities in \cref{thm:md}.
\end{remark}

The application of \cref{thm:md} to bandits and partial monitoring requires a few more ideas. The learner must not only choose the potential and learning rate(s), but also a way
of estimating the losses.  The latter is generally not possible without randomisation, so the learner must also choose a distribution from which to sample its actions.
To emphasise the presence of randomness, we now use capitals $(Q_t)_{t=1}^n$ for the recommendations of MD/FTRL and let $P_t \in \sP$ be the distribution from which the learner samples
action $A_t$. Very often $P_t$ has mean $Q_t$ but this is not universally true.
For example, in linear bandits $P_t$ is obtained by mixing $Q_t$ with a distribution on the contact points of John's ellipsoid or a Kiefer--Wolfowitz distribution \citep{BCK12}.
The generic outline of MD/FTRL as applied to bandits is given in \cref{alg:md}.

\begin{algorithm}
\textbf{input: } learning rate $\eta > 0$ and Legendre potential $F$ \\[-0.3cm]

\textbf{initialisation: } $Q_1 = \argmin_{q \in \cD} F(q)$ \\[-0.3cm]

\textbf{in each round $t$:} \\[-0.3cm]

$\qquad$\textbf{optimise:} compute exploratory distribution $P_t \in \sP$ based on history \\[-0.3cm]

$\qquad$\textbf{act:} sample $A_t \sim P_t$ and observe $\sigma_t = \Phi_{A_t}(z_t)$ \\[-0.3cm]

$\qquad$\textbf{update:} compute loss estimate $\hat \ell_t \in \R^d$ based on observations
\begin{align*}
\tag{MD} Q_{t+1} &= \argmin_{q \in \cD} \ip{q, \hat \ell_t} + \frac{1}{\eta} \KL(q, Q_t)  \\
\tag{FTRL} Q_{t+1} &= \argmin_{q \in \cD} \sum_{s=1}^t \ip{q, \hat \ell_s} + \frac{F(q)}{\eta} 
\end{align*}
\caption{Online stochastic MD/FTRL. Generally speaking, $P_t$ only depends on the history via $Q_t$. The learning rate is constant in the above, 
while in \cref{sec:adaptive} we will use an adaptive learning rates.}
\label{alg:md}
\end{algorithm}

\section{A generalised information ratio}

The information ratio was introduced by \cite{RV14} as a tool for the analysis of an algorithm called information-directed sampling, which explicitly optimises
the exploration/exploitation dilemma in a Bayesian framework. 
This beautiful idea led to a number of short proofs bounding the Bayesian regret for a variety of set-ups \citep{RV14,BDKP15,RV16,DV18,DMV19,LS19pminfo,Lat20-cvx}.
We introduce a generalisation of the concept and explore the properties of information-directed sampling.

\begin{definition}\label{def:info}
A partial monitoring game has a (generalised) information ratio of $(\alpha, \beta, \lambda)$ if 
for any $\nu \in \cV$, there exists a distribution $p \in \sP$ such that when $(Z, A^*, A)$ has law $\nu \otimes p$, then 
\begin{align*}
\E[\ip{A - A^*, \ell(Z)}] \leq \alpha + \beta^{1-1/\lambda}\E[\KL(\E[A^*|\Phi_A(Z), A], \E[A^*])]^{1/\lambda}\,.
\end{align*}
\end{definition}

The distributions $p \in \sP$ realising the display are called exploratory distributions.
The innovation in the following theorem is that previous work only addressed the case where $\lambda = 2$.

\begin{theorem}\label{thm:info}
Suppose a partial monitoring game has an information ratio of $(\alpha, \beta, \lambda)$ with $\alpha, \beta \geq 0$ and $\lambda \geq 1$. Then, for any finitely 
supported distribution $\mu$ on $\cZ^n$, there exists a policy $\pi$ such that
\begin{align*}
\BReg_n(\pi, \mu) \leq n (\epsd + \alpha) + (\beta n)^{1 -1/\lambda} \diam(\cD)^{1/\lambda}\,.
\end{align*}
\end{theorem}

The assumption that the prior $\mu$ is finitely supported is needed because in \cref{def:info} we only assumed the existence of a good exploratory distribution 
for distributions $\nu \in \cV$. 
Those concerned mostly with the Bayesian setting usually define the information ratio for a richer class of distributions than $\cV$ and correspondingly \cref{thm:info} would apply to more
priors. The reason for the choices here is for the connection to the stability term in \cref{thm:md}, where (a) the coarse $\cV$ is sufficient and (b) richer classes cause measure-theoretic challenges.

\begin{proof}[Proof of \cref{thm:info}]
Let $(Z_t)_{t=1}^n$ be the sequence of outcomes sampled from the prior $\mu$ and
$\E_t[\cdot]$ be the conditional expectation given the observation history $(A_s)_{s=1}^t$, $(\sigma_s)_{s=1}^t$ and abbreviate $\ell_t = \ell(Z_t)$.
Let
\begin{align*}
A^* = \argmin_{a \in \cD} \sum_{t=1}^n \ip{a, \ell_t}\,.
\end{align*}
Let $A^*_t = \E_{t-1}[A^*]$ be the expectation of the optimal action given the information available at the start of round $t$. Consider the policy $\pi$ that samples $A_t$ from any distribution $P_t \in \sP$ for which
\begin{align}
\E_{t-1}[\ip{A_t - A^*, \ell_t}] \leq \alpha + \beta^{1-1/\lambda} \E_{t-1}[\KL(A^*_{t+1}, A^*_t)]^{1/\lambda}\,,
\label{eq:inf-proof}
\end{align}
the existence of which is guaranteed by the assumptions of the theorem. 
Note, that here we have used the fact that $A_t$ and $(Z_t, A^*)$ are conditionally independent given $(A_s)_{s=1}^{t-1}$ and $(\sigma_s)_{s=1}^{t-1}$.
The Bayesian regret of this policy is bounded by
\begin{align*}
\BReg_n(\pi, \mu)
&= \E\left[\max_{a \in \cA} \sum_{t=1}^n \ip{A_t - a, \ell_t}\right] \\
&\leq n\epsd + \E\left[\sum_{t=1}^n \ip{A_t - A^*, \ell_t}\right] \\
&\leq n (\epsd + \alpha) + \E\left[\sum_{t=1}^n \beta^{1-1/\lambda} \E_{t-1}\left[\KL(A^*_{t+1}, A^*_t)\right]^{1/\lambda}\right] \\
&\leq n (\epsd + \alpha) + (\beta n)^{1 - 1/\lambda} \E\left[\sum_{t=1}^n \KL(A^*_{t+1}, A^*_t)\right]^{1/\lambda} \\
&\leq n (\epsd + \alpha) + (\beta n)^{1 - 1/\lambda} \diam(\cD)^{1/\lambda}\,,
\end{align*}
where the second inequality follows from \cref{eq:inf-proof}, the third from Jensen's inequality and the concavity of $x \mapsto x^{1/\lambda}$. The fourth inequality follows 
by telescoping the Bregman divergences \citep[Theorem 3]{LS19pminfo}.
\end{proof}

\paragraph{Information-directed sampling} 
Before moving on, let us explain the name `information ratio' and explore some properties of the information-directed sampling algorithm introduced by \cite{RV14} in 
the context of our generalisation.
Suppose that $(Z_t)_{t=1}^n$ are sampled from known finitely supported prior $\mu$ on $\cZ^n$ and 
$A^* = \argmin_{a \in \cD} \sum_{t=1}^n \ip{a, \ell(Z_t)}$.
Information-directed sampling is a Bayesian algorithm. In each round it solves an optimisation problem 
to find an exploratory distribution that minimises the ratio of the expected instantaneous squared regret
and the information gain, with the latter measured by the expected Bregman divergence between posterior and prior. 
This ratio is called the information ratio and the algorithm is summarised in \cref{alg:ids}.

For partial monitoring games with an information ratio of $(\alpha, \beta, \lambda)$ with $\alpha = 0$ and $\lambda = 2$, 
the information-directed sampling algorithm chooses exactly the exploratory distribution used in the proof of \cref{thm:info}, and
hence recovers the same bound. In light of the generalised definitions, however, one might question the extent to which
the optimisation problem in \cref{alg:ids} is fundamental. When $\lambda = 3$, the information ratio should perhaps be defined
as the cube of the regret divided by the information gain.
The next theorem provides an upper bound on the Bayesian regret of information-directed sampling that nearly matches \cref{thm:info}
for $\alpha = 0$ and $\lambda \geq 2$ without modifying the algorithm.

\begin{algorithm}
\textbf{input: } prior $\mu$ on $\cZ^n$

\textbf{for $t = 1$ to $n$:} \\[-0.3cm]

\qquad let $\E_{t-1}[\cdot] \triangleq \E[\cdot | A_1,\sigma_1,\ldots,A_{t-1},\sigma_{t-1}]$ \\[-0.2cm]

\qquad compute expected regret and information vectors:
\begin{align*}
\Delta_{t,a} &= \E_{t-1}[\ip{a - A^*, \ell(Z_t)}] \\
\cI_{t,a} &= \E_{t-1}[\KL(\E_{t-1}[A^*|\Phi_a(Z_t)], \E_{t-1}[A^*])]
\end{align*}
\qquad compute exploratory distribution:
\begin{align}
P_t = \argmin_{p \in \sP} \lim_{\epsilon \to 0^+} \frac{\max(0, \ip{p, \Delta_t})^2}{\ip{p, \cI_t} + \epsilon}
\label{eq:ids}
\end{align}
\qquad sample $A_t \sim P_t$ and observe $\sigma_t = \Phi_{A_t}(Z_t)$
\caption{Information-directed sampling}\label{alg:ids}
\end{algorithm}

\begin{theorem}\label{thm:ids}
Suppose a partial monitoring game has information ratio $(\alpha, \beta, \lambda)$ with $\alpha = 0$ and $\lambda \geq 2$.
Then the Bayesian regret of information-directed sampling is bounded for any finitely supported prior distribution $\mu$ on $\cZ^n$ by
\begin{align*}
\BReg_n \leq n \epsd + \lambda^{1-2/\lambda}  (\beta n)^{1-1/\lambda} \diam(\cD)^{1/\lambda}\,.
\end{align*}
\end{theorem}

Notice that (a) the theorem only holds for $\alpha = 0$ and $\lambda \geq 2$, 
(b) the algorithm does not depend on $\lambda$, 
and (c) the leading constant in \cref{thm:ids} is slightly worse than \cref{thm:info}.
The improved constant can be recovered by changing the optimisation problem in the definition of the algorithm to
\begin{align*}
P_t = \argmin_{p \in \sP} \lim_{\epsilon \to 0^+} \frac{\max(0, \ip{p, \Delta_t})^\lambda}{\ip{p, \cI_t} + \epsilon}\,.
\end{align*}
On the other hand, the resulting algorithm now depends on $\lambda$ and when $\lambda < 2$, the optimisation is not in general convex.

\begin{proof}[Proof of \cref{thm:ids}]
By the definition of the algorithm,
\begin{align*}
P_t = \argmin_{p \in \sP} \lim_{\epsilon \to 0^+} \frac{\max(0, \ip{p, \Delta_t})^2}{\epsilon + \ip{p, \cI_t}}\,.
\end{align*}
Suppose for a moment that $\ip{P_t, \Delta_t} > 0$. Then by the definition of the information ratio and \cref{lem:ids} in the appendix,
\begin{align*}
\frac{\ip{P_t, \Delta_t}^\lambda}{\ip{P_t, \cI_t}}
\leq 2^{\lambda - 2} \min_{p \in \sP} \frac{\ip{p, \Delta_t}^\lambda}{\ip{p, \cI_t}} \leq 2^{\lambda - 2} \beta^{\lambda - 1}\,.
\end{align*}
Therefore $\ip{P_t, \Delta_t} \leq 2^{1 - 2/\lambda} \beta^{1-1/\lambda} \ip{P_t, \cI_t}^\lambda$, which is obvious when $\ip{P_t, \Delta_t} \leq 0$. The Bayesian regret is now bounded using the same argument as in the proof of \cref{thm:info}:
\begin{align*}
\BReg_n 
&\leq n\epsd + \E\left[\sum_{t=1}^n \ip{P_t, \Delta_t}\right] \\
&\leq n\epsd + \E\left[\sum_{t=1}^n 2^{1-2/\lambda} \beta^{1-1/\lambda} \ip{P_t, \cI_t}^\lambda\right] \\
&\leq n\epsd + 2^{1-2/\lambda} (n \beta)^{1-1/\lambda} \E\left[\sum_{t=1}^n \ip{P_t, \cI_t}\right]^{1/\lambda} \\
&\leq n\epsd + 2^{1-2/\lambda} \diam(\cD)^{1/\lambda} (n \beta)^{1-1/\lambda}\,.
\qedhere
\end{align*}
\end{proof}

\begin{remark}
The variant of information-directed sampling using the Bregman divergence was introduced briefly by \cite{LS19pminfo}, generalising 
the original by \cite{RV14}, who used the mutual information.
The observation that information-directed sampling with a squared regret in the information ratio is reasonable even when $\lambda = 3$ was noticed already in the context of globally observable linear partial monitoring games by \cite{KLK20}.
\end{remark}

\section{Exploration by optimisation}\label{sec:exp-opt}

The policy introduced in this section uses the skeleton of \cref{alg:md} and solves an optimisation problem to find exploratory distributions and loss estimators in a way
that essentially minimises the bound. A similar algorithm has been seen before with a less clean form and in the context of finite partial monitoring \citep{LS19pmsimple}.

\paragraph{Optimisation problem}
Let $\cG$ be the space of functions from $\cA \times \Sigma$ to $\R^d$. Functions in $\cG$ will be used to estimate the losses and are called estimation functions.
An estimation function $g \in \cG$ is called unbiased if for all $z \in \cZ$ and $b, c \in \cA$,
\begin{align*}
\bip{b - c, \ell(z) - \sum_{a \in \cA} g(a, \Phi_a(z))} = 0\,.
\end{align*}
An unbiased loss estimation function $g \in \cG$ can be combined with importance-weighting to estimate relative differences in losses. 
Specifically, given any $p \in \sPp$ and $A \sim p$. Then, for any $z \in \cZ$,
\begin{align*}
\E\Bigg[\bip{b - c, \underbracket{\frac{g(A, \Phi_A(z))}{p(A)}}_{\smash[b]{\text{loss estimate}}}}\Bigg] = \ip{b - c, \ell(z)}\,.
\end{align*}
We now define the objective for an optimisation problem that plays a central role in everything that follows.
Given $q \in \cD \cap \dom(\nabla F)$ and $\eta > 0$, define a function $\Lambda_{q,\eta} : \cZ \times \cD \times \sPp \times \cG \to \R$ by
\begin{align*}
\Lambda_{q,\eta}(z, a^*, p, g) 
&= \sum_{a \in \cA} p(a) \bip{a - a^*, \ell(z)} + \bip{a^* - q, \sum_{a \in \cA} g(a, \Phi_a(z))} \\ &\quad\qquad + \frac{1}{\eta} \sum_{a \in \cA} p(a) \Psi_q\left(\frac{\eta g(a, \Phi_a(z))}{p(a)}\right) \,.
\end{align*}
Since sums of convex functions are convex and the perspective of a convex function is convex, the function $(p, g) \mapsto \Lambda_{q,\eta}(z, a^*, p,g)$ is convex.
To give a little more intuition for $\Lambda_{q,\eta}$, notice that
\begin{align*}
\Lambda_{q,\eta}(z, a^*, p, g)
&= \sum_{a \in \cA} p(a) \bip{a - q, \ell(z)} + \bip{a^* - q, \sum_{a \in \cA} g(a, \Phi_a(z)) - \ell(z)} \\ &\qquad \quad
+ \frac{1}{\eta} \sum_{a \in \cA} p(a) \Psi_q\left(\frac{\eta g(a, \Phi_a(z))}{p(a)}\right) \,.
\end{align*}
The first term measures the loss due to sampling an action from $p$ with mean $\sum_{a \in \cA} p(a) a$ rather than a distribution with mean $q$ as recommended by MD/FTRL.
The second term vanishes when $g$ is unbiased and otherwise provides some measure of the bias.
The last term measures the stability of the online learning algorithm.
Define $\Lambda^*_{q,\eta}$ and $\Lambda^*_\eta$ by
\begin{align}
\Lambda^*_{q,\eta} &= \inf_{\substack{p \in \sPp \\ g \in \cG}} \sup_{\substack{z \in \cZ \\ a^* \in \cD}} \Lambda_{q,\eta}(z, a^*, p, g) &
\Lambda^*_\eta &= \sup_{q \in \cD \cap \dom(\nabla F)} \Lambda^*_{q,\eta} \,. \label{eq:Lambda}
\end{align}

\begin{algorithm}
\textbf{input: } Learning rate $\eta$ and precision $\epsilon$

\noindent \textbf{initialise: } $Q_1 = \argmin_{q \in \cD} F(q)$ 

\noindent \textbf{for $\bm{t = 1}$ to $\bm{n}$:} \\[-0.2cm]

$\qquad$\textbf{optimisation:} find exploratory distribution $P_t \in \sPp$ and $G_t \in \cG$ such that
\begin{align*}
\sup_{z \in \cZ, a^* \in \cD} \Lambda_{Q_t,\eta}(z, a^*, P_t, G_t) \leq \Lambda^*_\eta + \epsilon 
\end{align*}

$\qquad$\textbf{acting:}  sample action $A_t \sim P_t$ and observe signal $\sigma_t = \Phi_{A_t}(z_t)$ \\

$\qquad$\textbf{update:} compute loss estimate and $Q_{t+1}$
\begin{align*}
\hat \ell_t &= \frac{G_t(A_t, \sigma_t)}{P_t(A_t)} \\
\tag{MD} Q_{t+1} &= \argmin_{q \in \cD} \ip{q, \hat \ell_t} + \frac{1}{\eta} \KL(q, Q_t) \\
\tag{FTRL} Q_{t+1} &= \argmin_{q \in \cD} \sum_{s=1}^t \ip{q, \hat \ell_s} + \frac{F(q)}{\eta} 
\end{align*}
\caption{Exploration by optimisation}\label{alg:exp-opt}
\end{algorithm}

\begin{theorem}\label{thm:exp-opt}
The regret of the policy defined by \cref{alg:exp-opt} (using either MD or FTRL) when run with precision $\epsilon > 0$ and learning rate $\eta > 0$ is bounded by
\begin{align*}
\Reg_n \leq \frac{\diam(\cD)}{\eta} + n(\epsd + \epsilon + \Lambda^*_\eta) \,.
\end{align*}
\end{theorem}

\begin{proof}
Let $\ell_t = \ell(z_t)$ and $a^* = \argmin_{a \in \cD} \sum_{t=1}^n \ip{a, \ell_t}$ be the optimal action in hindsight. Decomposing the regret relative to $a^*$ and applying \cref{thm:md} yields 
\begin{align*}
\Reg_n
&\leq n\epsd + \E\left[\sum_{t=1}^n \sum_{a \in \cA} P_t(a) \ip{a - a^*, \ell_t}\right] \\
&= n\epsd + \E\left[\sum_{t=1}^n \sum_{a \in \cA} P_t(a) \ip{a - a^*, \ell_t} + \ip{a^* - Q_t, \hat \ell_t} + \ip{Q_t - a^*, \hat \ell_t}\right] \\
&\leq n\epsd + \frac{\diam(\cD)}{\eta} + \sum_{t=1}^n 
\underbracket{\E\left[\sum_{a \in \cA} P_t(a)\ip{a - a^*, \ell_t} + \ip{a^* - Q_t, \hat \ell_t} + \frac{1}{\eta} \Psi_{Q_t}(\eta \hat \ell_t)\right]}_{\textrm{(A)}_t}\,.
\end{align*}
Using the fact that $p \in \sPp$ and the definition of expectation yields 
\begin{align*}
\E[\textrm{(A)}_t] 
&= \E\Bigg[\sum_{a \in \cA} P_t(a) \ip{a - a^*, \ell_t} + \bip{a^* - Q_t, \sum_{a \in \cA} G_t(a, \Phi_a(z_t)) } \\
& \qquad \qquad+ \frac{1}{\eta} \sum_{a \in \cA} P_t(a) \Psi_{Q_t}\left(\frac{\eta G_t(a, \Phi_a(z_t))}{P_t(a)}\right) \Bigg] \\
&= \E\left[\Lambda_{Q_t,\eta}(z_t, a^*, P_t, G_t)\right] \\
&\leq \Lambda^*_{\eta} + \epsilon \,, 
\end{align*}
where the last inequality follows from the definition of $P_t$ and $G_t$ in \cref{alg:exp-opt}.
\end{proof}

\section{Stability and the information ratio}\label{sec:inf-stab}

The next theorem makes a connection between the information ratio and the value of the optimisation problems defined in \cref{eq:Lambda}.

\begin{theorem}\label{thm:inf}
Suppose a partial monitoring game has an information ratio of $(\alpha, \beta, \lambda)$ with $\lambda > 1$. Then,
\begin{align*}
\Lambda^*_\eta \leq \alpha + \beta \left(1 - \frac{1}{\lambda}\right) \left(\frac{\eta}{\lambda}\right)^{\frac{1}{\lambda - 1}} \,.
\end{align*}
\end{theorem}

\begin{corollary}
The regret of \cref{alg:exp-opt} with precision $\epsilon > 0$ and learning rate
\begin{align*}
\eta = \lambda \left(\frac{\diam(\cD)}{\beta n}\right)^{1 - 1/\lambda}
\end{align*}
is bounded by $\displaystyle \Reg_n \leq (\epsilon + \epsd + \alpha) n +  \diam(\cD)^{\frac{1}{\lambda}} (\beta n)^{1 - \frac{1}{\lambda}}$.
\end{corollary}

\begin{proof}
Combine \cref{thm:exp-opt,thm:inf}.
\end{proof}

Before the proof of \cref{thm:inf}, we start with a technical lemma lower bounding $\Lambda_{q,\eta}$.

\begin{lemma}\label{lem:bound}
Let $q \in \cD \cap \dom(\nabla F)$ and $\eta > 0$. Then there exists a constant $C_{q,\eta}$ such that
\begin{align*}
\Lambda_{q,\eta}(z, a^*, p, g) \geq C_{q,\eta} 
\end{align*}
for all $z \in \cZ$, $a^* \in \cD$, $p \in \sPp$ and $g \in \cG$.
\end{lemma}

\begin{proof}
By the Fenchel--Young inequality,
\begin{align*}
&\bip{a^* - q, \frac{g(a, \Phi_a(z))}{p(a)}} + \frac{1}{\eta} \Psi_q\left(\frac{\eta g(a, \Phi_a(z))}{p(a)}\right) \\ 
&= \frac{1}{\eta} \bip{a^*, \frac{\eta g(a, \Phi_a(z))}{p(a)}} + \frac{1}{\eta} F^*\left(\nabla F(q) - \frac{\eta g(a, \Phi_a(z))}{p(a)}\right) - \frac{1}{\eta} F^*(\nabla F(q)) \\
&\geq \frac{\ip{a^*, \nabla F(q)} - F(a^*) - F^*(\nabla F(q))}{\eta} \\
&\geq -\frac{\norm{a^*} \norm{\nabla F(q)} + F(a^*) + F^*(\nabla F(q))}{\eta}\,.
\end{align*}
Hence, using the definition of $\Lambda_{q,\eta}$,
\begin{align*}
\Lambda_{q,\eta}(z, a^*, p, g) 
&\geq \sum_{a \in \cA} p(a) \ip{a - a^*, \ell(z)} - \frac{\norm{a^*} \norm{\nabla F(q)} + F(a^*) + F^*(\nabla F(q))}{\eta} \\
&\geq -\frac{\norm{a^*} \norm{\nabla F(q)} + F(a^*) + F^*(\nabla F(q))}{\eta} - 1\,,
\end{align*}
where in the last inequality we used the assumption that the losses are in $[0,1]$.
The right-hand side is lower bounded by a constant that depends only on $q$ and $\eta$ since $a^* \in \cD$ and $F$ has finite diameter on $\cD$.
\end{proof}

\begin{proof}[Proof of \cref{thm:inf}]
The core ingredients of the proof are an application of Sion's minimax theorem to exchange the $\inf$ and $\sup$ in the definition of $\Lambda^*_{q,\eta}$ and 
an algebraic calculation to introduce the information ratio.  
The argument is complicated by the fact that $\nabla F$ need not exist on $\cD \setminus \relint(\cD)$. 

\paragraph{Step 1: Notation and setup}
Let $q \in \cD \cap \dom(\nabla F)$ and $\eta > 0$ be fixed and abbreviate $\Lambda(\cdot) \equiv \Lambda_{q,\eta}(\cdot)$. 
Let $y \in \relint(\cD)$, which exists by assumption. For $\epsilon > 0$ define 
\begin{align*}
\cD_\epsilon = \{(1 - \epsilon) x + \epsilon y : x \in \cD\} \subset \relint(\cD)\,.
\end{align*}
Convexity of $F$ and the assumption that $\cD \subset \dom(F)$ and that $F$ is Legendre implies that $\sup_{x \in \cD_\epsilon} \norm{\nabla F(x)}_\infty < \infty$.
Let $\cV_\epsilon \subset \cV$ be the space of finitely supported probability distributions on $\cZ \times \cD_\epsilon$
and $\cG_\epsilon \subset \cG$ be the set of estimation functions $g$ with $\max_{a \in \cA} \sup_{\sigma \in \Sigma} \norm{g(a, \sigma)}_\infty \leq C_\epsilon$ where 
\begin{align*}
C_\epsilon = \frac{1}{\eta}\sup_{q' \in \cD_\epsilon} \norm{\nabla F(q) - \nabla F(q')}_\infty \,.
\end{align*}
Next, let $\cX_\epsilon \subset \sP_\epsilon \times \cG_\epsilon$ be given by
\begin{align*}
\cX_\epsilon = \left\{(p, g) : \nabla F(q) - \frac{\eta g(a, \sigma)}{p(a)} \in \nabla F(\cD_\epsilon) \text{ for all } \sigma \in \Sigma \right\}\,,
\end{align*}
which is convex (since $F$ is continuously differentiable by our assumptions) and compact.

\paragraph{Step 2: Exchanging inf and sup}
We will now use Sion's theorem to exchange the $\inf$ and the $\sup$ in the definition of $\Lambda$ and show that
\begin{align}
\inf_{\substack{p \in \sP \\ g \in \cG}} \sup_{\substack{a^* \in \cD \\ z \in \cZ}} \Lambda(z, a^*, p, g)
\leq \liminf_{\epsilon \to 0} \sup_{\nu \in \cV_\epsilon} \inf_{(p, g) \in \cX_\epsilon} \int_{\cZ \times \cD} \Lambda(z, a^*, p, g) \d{\nu}(z, a^*)\,.
\label{eq:sion}
\end{align}
The analysis in this step depends on some topological tomfoolery and can be skipped by eager readers.
Imbue $\cG_\epsilon$ with the product topology, which is the initial topology of the collection of maps $(g \mapsto g(a, \sigma))_{a \in \cA, \sigma \in \Sigma}$.
In other words, the topology on $\cG_\epsilon$ is the coarsest topology such that $g \mapsto g(a, \sigma)$ is continuous for all $a \in \cA$ and $\sigma \in \Sigma$.
By Tychonoff's theorem, $\cG_\epsilon$ is compact while $\sP_\epsilon$ is compact with the usual topology.
Furthermore, when $(p, g) \in \cX_\epsilon$,
\begin{align*}
\sup_{z \in \cZ, a^* \in \cD_\epsilon} \Lambda(z, a^*, p, g) < \infty\,.
\end{align*}
In combination with \cref{lem:bound}, this shows that $\Lambda$ is bounded on the domain $\cZ \times \cD_\epsilon \times \cX_\epsilon$.
Continuity of $(p, g) \mapsto \Lambda(a^*, z, p, g)$ follows from the definition of the product topology and the same mapping is convex via the perspective construction as noted in \cref{sec:exp-opt}.
By choosing the discrete topology on $\cZ \times \cD_\epsilon$, the mapping
$(z, a^*) \mapsto \Lambda(a^*, z, p, g)$ is automatically continuous. 
Let $\cV_\epsilon$ have the weak* topology and $(p, g) \in \cX_\epsilon$.
Then $\nu \mapsto \int_{\cZ \times \cD_\epsilon} \Lambda(z, a^*, p, g) \d{\nu}(z, a^*)$ is continuous by the definition of the weak* topology and using the previous argument that $\Lambda(\cdot, \cdot, p, g)$ is bounded when $(p, g) \in \cX_\epsilon$. The same mapping is clearly linear.
Hence, by Sion's minimax theorem \citep{Sio58}, 
\begin{align*}
\inf_{(p, g) \in \cX_\epsilon} \sup_{\substack{a^* \in \cD_\epsilon \\ z \in \cZ}}  \Lambda(z, a^*, p, g)
&= \inf_{(p, g) \in \cX_\epsilon} \sup_{\nu \in \cV_\epsilon} \int_{\cZ \times \cD_\epsilon} \Lambda(z, a^*, p, g) \d{\nu}(z, a^*) \\
&= \sup_{\nu \in \cV_\epsilon} \inf_{(p, g) \in \cX_\epsilon} \int_{\cZ \times \cD_\epsilon} \Lambda(z, a^*, p, g) \d{\nu}(z, a^*) \,.
\end{align*}
Combining this with linearity of the map $a \mapsto \Lambda(z, a, p, g)$ and \cref{lem:bound} shows that
\begin{align*}
\inf_{\substack{p \in \sP \\ g \in \cG}} \sup_{\substack{a^* \in \cD \\ z \in \cZ}} \Lambda(z, a^*, p, g)
&\leq \inf_{(p, g) \in \cX_\epsilon} \sup_{\substack{a^* \in \cD \\ z \in \cZ}} \Lambda(z, a^*, p, g) \\
&= \inf_{(p, g) \in \cX_\epsilon} \sup_{\substack{a^* \in \cD \\ z \in \cZ}} \frac{\Lambda(z, \epsilon y + (1 - \epsilon) a^*, p, g) - \epsilon \Lambda(z, y, p, g)}{1 - \epsilon} \\
&\leq \inf_{(p, g) \in \cX_\epsilon} \sup_{\substack{a^* \in \cD_\epsilon \\ z \in \cZ}} \frac{\Lambda(z, a^*, p, g) - \epsilon C_{q,\eta}}{1 - \epsilon} \\
&= \sup_{\nu \in \cV_\epsilon} \inf_{(p, g) \in \cX_\epsilon} \int_{Z \times \cD} \frac{\Lambda(z, a^*, p, g)}{1 - \epsilon} \d{\nu}(z, a^*) - \frac{\epsilon C_{q, \eta}}{1 - \epsilon}
\end{align*}
Taking the limit as $\epsilon$ tends to zero establishes \cref{eq:sion}.

\paragraph{Step 3: Introducing the information ratio}
Fix $\epsilon > 0$ and $\nu \in \cV_\epsilon$ and $p \in \sP_\epsilon$ and let $(Z,A^*, A)$ have law $\nu \otimes p$, which means that
\begin{align*}
&\int_{Z \times \cD} \Lambda(z, a^*, p, g) \d{\nu}(z, a^*)
= \E[\Lambda(Z, A^*, p, g)] \\
&\qquad= \E\left[\bip{A - A^*, \ell(Z)} + \bip{A^* - q, \frac{g(A, \sigma)}{p(A)}} + \frac{1}{\eta} \Psi_q\left(\frac{\eta g(A, \sigma)}{p(A)}\right)\right]\,.
\end{align*}
The first term will be bounded using \cref{lem:supp} and the assumptions on the information ratio.
The second term is bounded by explicitly minimising the second term.
Given any action $a \in \cA$ and signal $\sigma \in \{\Phi_a(z) : z \in \cZ,\, \nu(\{z\} \times \cD_\epsilon) > 0\}$, let
\begin{align*}
g(a, \sigma) = \frac{p(a)}{\eta} \left(\nabla F(q) - \nabla F(\E[A^* | \Phi_a(Z) = \sigma])\right)\,,
\end{align*}
and otherwise let $g(a, \sigma) = 0$. 
Since $A^* \in \cD_\epsilon$, it holds that $\E[A^* | \Phi_a(Z) = \sigma] \in \cD_\epsilon$. Therefore
$\max_{a \in \cA}, \sup_{\sigma \in \Sigma} \norm{g(a, \sigma)}_\infty \leq C_\epsilon$, which implies that $g \in \cG_\epsilon$ and hence $(p, g) \in \cX_\epsilon$.
Next, let $\sigma = \Phi_A(Z)$ and $\prior = \E[A^*]$ and $\post = \E[A^* | A, \Phi_A(Z)]$. 
Then, using the definitions, non-negativity of the Bregman divergences and duality (\cref{eq:kl-dual}), 
\begin{align}
&\E\left[\bip{A^* - q, \frac{g(A, \sigma)}{p(A)}} + \frac{1}{\eta} \Psi_q\left(\frac{\eta g(A, \sigma)}{p(A)}\right)\right] \nonumber \\
&\qquad= \frac{1}{\eta} \E\left[\ip{A^*, \nabla F(q) - \nabla F(\post)} + F^\star\left(\nabla F(\post)\right) - F^\star(\nabla F(q))\right] \nonumber \\
&\qquad= -\frac{1}{\eta}\E\left[F^\star(\nabla F(q)) - F^\star(\nabla F(\prior)) - \ip{A^*, \nabla F(q) - \nabla F(\prior)}\right] \nonumber \\
&\qquad\qquad -\frac{1}{\eta}\E\left[F^\star(\nabla F(\prior)) - F^\star(\nabla F(\post)) - \ip{A^*, \nabla F(\prior) - \nabla F(\post)}\right] \nonumber \\
&\qquad= -\frac{1}{\eta}\E\left[F^\star(\nabla F(q)) - F^\star(\nabla F(\prior)) - \ip{\prior, \nabla F(q) - \nabla F(\prior)}\right] \nonumber \\
&\qquad\qquad -\frac{1}{\eta}\E\left[F^\star(\nabla F(\prior)) - F^\star(\nabla F(\post)) - \ip{\post, \nabla F(\prior) - \nabla F(\post)}\right] \nonumber \\
&\qquad= -\frac{1}{\eta}\E\left[\KL_\star(\nabla F(q), \nabla F(\prior)) + \KL_\star(\nabla F(\prior), \nabla F(\post))\right] \nonumber \\
&\qquad\leq -\frac{1}{\eta}\E\left[\KL_\star(\nabla F(\prior), \nabla F(\post))\right] \nonumber \\
&\qquad= -\frac{1}{\eta}\E\left[\KL(\post, \prior)\right]\,, \label{eq:dual}
\end{align}
where the first equality follows from the definitions of $\Psi_q$ and $g$. Note, $g$ was chosen so as to minimise this expression.
The second by adding and subtracting terms.
The third is true by the definition of $\prior$ and $\post$ and the fourth is the definition of the Bregman divergence.
The inequality is true since Bregman divergences are always non-negative. The final equality follows from duality (\cref{eq:kl-dual}).
By \cref{lem:supp}, $p \in \sP_\epsilon$ can be chosen so that 
\begin{align}
\E\left[\ip{A - A^*, \ell(Z)}\right] \leq \epsilon + \alpha + \beta^{1-1/\lambda} \E\left[\KL(\post,\prior)\right]^{1/\lambda} \,.
\label{eq:dist}
\end{align}
Combining this with \cref{eq:dual}, the definition of $\Lambda$ and elementary optimisation shows that
\begin{align*}
\E\left[\Lambda(Z, A^*, p, g)\right] 
&\leq \E\left[\ip{A - A^*, \ell(Z)} - \frac{1}{\eta} \KL(\post, \prior)\right] \\
&\leq |\cA| \epsilon + \alpha + \beta^{1-1/\lambda} \E[\KL(\post, \prior)]^{1/\lambda} - \frac{1}{\eta} \E[\KL(\post,\prior)] \\
&\leq |\cA| \epsilon + \alpha + \beta \left(1 - \frac{1}{\lambda}\right) \left(\frac{\eta}{\lambda}\right)^{\frac{1}{\lambda - 1}}\,.
\end{align*}
All together we have shown that for any $\epsilon > 0$ and $\nu \in \cV_\epsilon$ there exists a $(p, g) \in \cX_\epsilon$ such that
\begin{align*}
\int_{\cZ \times \cD} \Lambda(z, a^*, p, g) \d{\nu}(z, a^*) 
\leq |\cA| \epsilon + \alpha + \beta \left(1 - \frac{1}{\lambda}\right) \left(\frac{\eta}{\lambda}\right)^{\frac{1}{\lambda - 1}}\,.
\end{align*}
The claim of the theorem now follows from \cref{eq:sion}.
\end{proof}

\cref{thm:inf} provides a bound on $\Lambda^*_\eta$ in terms of the information ratio, but does not provide much information about which policy and estimation
functions that yield the bound. 
A fundamental case where more information can be extracted is when $\cA = \{e_1,\ldots,e_d\}$ and a bound on the information ratio
is witnessed by Thompson sampling, as is often the case. 
The next theorem relies on a class of potential functions that are widely used in finite-armed bandits \citep[for example]{CL18,ZiSe19}. 
Given $s \in \R$, the $s$-Tsallis entropy is
\begin{align*}
F(p) = \sum_{a=1}^d \frac{p_a^s - s p_a - (1-s)}{s(s-1)}\,.
\end{align*}
The limits as $s \to 1$ and $s \to 0$ correspond to the negentropy and logarithmic barrier, respectively.

\begin{theorem}\label{thm:ts}
Suppose that $F$ is the $s$-Tsallis entropy with $s \in [0,1]$ and 
$\cA = \{e_1,\ldots,e_d\}$ and $\cD = \sP_{\eta^{4/3}}$. Assume that for any $(Z, A^*)$ with law $\nu \in \cV$ and independent $A$ with law $p = \E[A^*]$,
\begin{align*}
\E\left[\left|\E[\ip{A, \ell(Z)}] - \E[\ip{A^*, \ell(Z)}]\right| \right]
\leq \sqrt{\beta \E[\KL(\E[A^*|\Phi_A(Z), A], \E[A^*])]}\,,
\end{align*}
where $\beta \geq 0$ is a constant. Then, 
\begin{align*}
\inf_{g \in \cG} \sup_{\substack{a^* \in \cD \\ z \in \cZ}} \Lambda_{q,\eta}(z, a^*, p, g) \leq (1 + O(\eta^{2/3})) \frac{\beta \eta}{4}\,,
\end{align*}
where the Big-O hides a constant depending only on $\beta$.
\end{theorem}

Note, the presence of the absolute values in the conditions of \cref{thm:ts} is slightly stronger than the definition of the information ratio in \cref{def:info}.
As far as we are aware, all known bounds on the information ratio hold for this stronger definition.

\begin{corollary}
Under the same assumptions as \cref{thm:ts},
there exist estimation functions such that MD/FTRL with $\cD = \sP_{\eta^{4/3}}$ and $P_t = Q_t$ and 
\begin{align*}
\eta = 2 \sqrt{\frac{\diam(\cD)}{n \beta}}
\end{align*}
satisfies $\Reg_n = \sqrt{(1 + o(1)) \beta n \diam(\cD)}$.
\end{corollary}

\begin{proof}
Combine \cref{thm:exp-opt,thm:ts} yields the following corollary and note that $\epsd \leq d \eta^{4/3}$, which contributes negligibly for large $n$. 
\end{proof}

\section{Adaptivity}\label{sec:adaptive}

Data-dependent analysis of bandit algorithms based on exponential weights or FTRL has a long history \citep[for example]{AAGO06}.
Recently, \cite{BS20} developed a data-dependent version of the information-theoretic analysis that was specified towards proving first-order bounds
for combinatorial semi-bandits. Here we generalise this concept by introducing an adaptive generalised information ratio and extending the results of earlier sections by showing
the existence of a corresponding FTRL strategy.

\begin{definition}
Let $\alpha \in \R$ and $\beta : \cZ \times \cA \to [0, \infty)$ and $\lambda > 1$. A partial monitoring game has an $(\alpha, \beta, \lambda)$ adaptive information ratio if
for all $\nu \in \cV$ there exists a $p \in \sP$ such that when $(Z, A^*, A)$ has law $\nu \otimes p$, then
\begin{align*}
\E[\ip{A - A^*, \ell(Z)}] \leq \alpha + \E[\beta(Z, A)]^{1-1/\lambda} \E[\KL(\E[A^*|\Phi_A(Z), A], \E[A^*])]^{1/\lambda}\,.
\end{align*}
\end{definition}

The next theorem is a straightforward generalisation of \cref{thm:info}. That theorem is recovered exactly when $\beta$ is a constant function.

\begin{theorem} \label{thm:adaptive-info}
Suppose a partial monitoring game has a $(\alpha, \beta, \lambda)$ adaptive information ratio, then for any prior $\nu \in \cV$, there exists a policy such that
\begin{align*}
\BReg_n \leq n(\epsd + \alpha) + \diam(\cD)^{1/\lambda} \E\left[\sum_{t=1}^n \beta(Z_t, A_t)\right]^{1-1/\lambda} \,,
\end{align*}
where $(Z_t)_{t=1}^n$ is sampled from $\nu$.
\end{theorem}

\begin{proof}
Using the same notation and argument as in \cref{thm:info}, 
\begin{align*}
\BReg_n 
&\leq n (\epsd + \alpha) + \E\left[\sum_{t=1}^n \E_{t-1}[\beta(Z_t, A_t)]^{1-1/\lambda} \E_{t-1}[\KL(A^*_{t+1}, A^*_t)]^{1/\lambda}\right] \\
&\leq n (\epsd + \alpha) + \diam(\cD)^{1/\lambda} \E\left[\sum_{t=1}^n \beta(Z_t, A_t)\right]^{1-1/\lambda} \,.
\qedhere
\end{align*}
\end{proof}

The next theorem generalises \cref{thm:inf}.

\begin{theorem} \label{thm:adaptive-inf}
Suppose a partial monitoring game has an $(\alpha, \beta, \lambda)$ adaptive information ratio and $\beta$ is bounded. Then, for any $\eta > 0$ and $q \in \cD \cap \dom(\nabla F)$,
\begin{align*}
\inf_{\substack{p \in \sPp \\ g \in \cG}} \sup_{\substack{a^* \in \cD \\ z \in \cZ}} \left[\Lambda_{q,\eta}(a^*, z, p, g) - \left(1 - \frac{1}{\lambda}\right) \left(\frac{\eta}{\lambda}\right)^{\frac{1}{\lambda - 1}} \sum_{a \in \cA} p(a) \beta(z, a)\right] \leq \alpha\,.
\end{align*}
\end{theorem}

\begin{proof}
Let $(Z, A^*, A)$ be the projection random element on measurable space $\cZ \times \cD \times \cA$ and $\E_{\nu,p}$ 
be the expectation with respect to probability measure $\nu \otimes p$ where $\nu \in \cV$ and $p \in \sPp$.
Given $\nu \in \cV$ and $p \in \sPp$, let 
\begin{align*}
\bar D_{\nu,p} &= \E_{\nu,p} [\KL(\E[A^* | \Phi_A(Z), A], \E[A^*])] &
\bar \beta_{\nu,p} &= \E_{\nu,p}[\beta(Z, A)] \,. 
\end{align*}
Notice that the term added inside the saddle point problem in the theorem statement is linear in $p$ and bounded by assumption.
Hence, the application of minimax theorem in the proof of \cref{thm:inf} goes through in the same manner, which shows that
\begin{align*}
&\inf_{p \in \sPp, g \in \cG} \sup_{a^* \in \cD, z \in \cZ} \Lambda_{q,\eta}(a^*, z, p, g) - 
  \left(1 - \frac{1}{\lambda}\right) \left(\frac{\eta}{\lambda}\right)^{\frac{1}{\lambda - 1}} \sum_{a \in \cA} p(a) \beta(z, a) \\ 
&\quad\leq \sup_{\nu \in \cV} \inf_{p \in \sPp} \left(\E_{\nu,p}\left[\ip{A - A^*, \ell(Z)}\right] - \left(1 - \frac{1}{\lambda}\right) \left(\frac{\eta}{\lambda}\right)^{\frac{1}{\lambda - 1}} \bar \beta_{\nu,p} - \frac{\bar D_{\nu,p}}{\eta}\right) \\
&\quad\leq \sup_{\nu \in \cV} \left(\alpha + \bar \beta_{\nu,p(\nu)}^{1-1/\lambda} \bar D_{\nu,p(\nu)}^{1/\lambda} - \left(1 - \frac{1}{\lambda}\right) \left(\frac{\eta}{\lambda}\right)^{\frac{1}{\lambda - 1}} \bar \beta_{\nu,p(\nu)} - \frac{\bar D_{\nu,p(\nu)}}{\eta}\right) \\
&\quad\leq \alpha\,,
\end{align*}
where the last inequality follows from elementary optimisation and
$p : \cV \to \sPp$ is a mapping guaranteed by the adaptive information ratio for which
\begin{align*}
\E_{\nu,p(\nu)}[\ip{A - A^*, \ell(Z)}] &\leq \alpha + \bar \beta_{\nu,p(\nu)}^{1-1/\lambda} \bar D_{\nu,p(\nu)}^{1/\lambda}\,. \qedhere
\end{align*}
\end{proof}

\cref{alg:exp-opt} can be made adaptive by optimising $P_t$ and $G_t$ so that
\begin{align*}
\sup_{a^* \in \cD, z \in \cZ} \Lambda_{Q_t,\eta}(a^*, z, P_t, G_t) - \left(1 - \frac{1}{\lambda}\right) \left(\frac{\eta }{\lambda}\right)^{\frac{1}{\lambda - 1}} \sum_{a \in \cA} P_t(a) \beta(z, a) \leq \alpha + \epsilon\,.
\end{align*}
By repeating the analysis in the proof of \cref{thm:exp-opt}, it follows that
\begin{align}
\Reg_n 
&\leq n (\epsilon + \epsd + \alpha) + \frac{\diam(\cD)}{\eta} + \left(1 - \frac{1}{\lambda}\right) \left(\frac{\eta}{\lambda}\right)^{\frac{1}{\lambda - 1}} \E\left[\sum_{t=1}^n \beta(z_t, A_t)\right]\,. 
\end{align}
There are two problems.
First, the expectation in the right-hand side depends on the law of the actions of the algorithm, which depend on $\eta$. Hence, it is not straightforward to optimise the learning rate.
Second, even if $(z, a) \mapsto \beta(z, a)$ can be written as a function of $z$ only, the quantity in the expectation is generally not known to the learner in advance. 
Both problems are resolved by tuning the learning rate online.

\paragraph{Online tuning}
Adaptively tuning the learning rate is possible if $(z, a) \mapsto \beta(z, a)$ can be written as a function of the signal $\Phi_a(z)$ and $a$.
For the remainder of the section we assume this is true and abuse notation by writing $\beta(\sigma, a)$.
Let
\begin{align}
\eta_t = \lambda^{-1/\lambda} (\lambda - 1)^{1-1/\lambda} \left(\frac{\diam(\cD)}{\beta_0 + \sum_{s=1}^{t-1} \beta(\sigma_s, A_s)} \right)^{1 - 1/\lambda} \,,
\label{eq:adaptive-eta}
\end{align}
where $\beta_0 = \sup_{\sigma \in \Sigma} \max_{a \in \cA} \beta(\sigma, a)$.
Consider the policy that chooses $P_t \in \sPp$ and $G_t \in \cG$ such that
\begin{align}
\sup_{\substack{z \in \cZ \\ a^* \in \cD}} \Lambda_{\eta_t,Q_t}(z, a^*, P_t, G_t) - 
\left(1 - \frac{1}{\lambda}\right) \left(\frac{\eta_t }{\lambda}\right)^{\frac{1}{\lambda - 1}} 
\sum_{a \in \cA} P_t(a) \beta(\Phi_a(z), a) \leq \epsilon + \alpha\,,
\label{eq:md-adaptive}
\end{align}
where $\eta_t$ is defined in \cref{eq:adaptive-eta} and with $\hat \ell_s = G_s(A_s, \sigma_s)$,
\begin{align*}
Q_t = \argmin_{q \in \cD} \sum_{s=1}^{t-1} \ip{q, \hat \ell_s} + \frac{F(q)}{\eta_t}\,.
\end{align*}
\begin{remark}
Mirror descent can behave badly when the learning rate is non-constant, so only the FTRL version of the algorithm is used here.
\end{remark}

\begin{theorem}\label{thm:adaptive-online}
The regret of the policy choosing $P_t$ and $G_t$ satisfying \cref{eq:md-adaptive} is bounded by
\begin{align*}
\Reg_n \leq n (\epsilon + \epsd + \alpha) + \left(\frac{\lambda}{\lambda-1}\right)^{1 - \frac{1}{\lambda}} \diam(\cD)^{\frac{1}{\lambda}} \E\left[\left(\beta_0 + \sum_{t=1}^{n-1} \beta(\sigma_t, A_t)\right)^{1-\frac{1}{\lambda}}\right]\,.
\end{align*}
\end{theorem}

\begin{proof}
Repeat the analysis in \cref{thm:exp-opt} to show that
\begin{align*}
\Reg_n \leq 
n (\epsilon + \epsd + \alpha) + \E\left[\frac{\diam(\cD)}{\eta_n} + \left(1 - \frac{1}{\lambda}\right) \sum_{t=1}^n \left(\frac{\eta_t}{\lambda}\right)^{1-1/\lambda} \beta(\sigma_t, A_t)\right]\,.
\end{align*}
Then combine the definition of $\eta_t$ with \cref{lem:tech1} in the appendix.
\end{proof}

The order of the expectation and $x \mapsto x^{1 - 1/\lambda}$ has been reversed in \cref{thm:adaptive-online} relative to \cref{thm:adaptive-info}, which except for the 
marginally larger leading constant and the presence of $\beta_0$ is actually an improvement. A similar improvement is possible in \cref{thm:adaptive-info}. Let $(\eta_t)_{t=1}^n$
be the sequence of learning rates as defined in \cref{eq:adaptive-eta}. Then, using the notation in the proof of \cref{thm:adaptive-info},
\begin{align*}
\BReg_n 
&\leq n (\epsd + \alpha) + \E\left[\sum_{t=1}^n \E_{t-1}[\beta(\sigma_t,A_t)]^{1-1/\lambda} \E_{t-1}[\KL(A^*_{t+1},A^*_t)]^{1/\lambda}\right] \\
&\leq n (\epsd + \alpha) + \E\left[\sum_{t=1}^n \frac{\E_{t-1}[\KL(A^*_{t+1}, A^*_t)]}{\eta_t} + (1 - \lambda) \left(\frac{\eta_t}{\lambda}\right)^{\frac{1}{\lambda - 1}} \beta(\sigma_t, A_t)\right] \\
&\leq n (\epsd + \alpha) + \E\left[\sum_{t=1}^n \frac{F(A^*_{t+1}) - F(A^*_t)}{\eta_t} + (1 - \lambda) \left(\frac{\eta_t}{\lambda}\right)^{\frac{1}{\lambda - 1}} \beta(\sigma_t, A_t)\right] \\
&\leq n (\epsd + \alpha) + \E\left[\frac{\diam(\cD)}{\eta_n} + (1 - \lambda) \sum_{t=1}^n \left(\frac{\eta_t}{\lambda}\right)^{\frac{1}{\lambda - 1}} \beta(\sigma_t, A_t)\right] \\
&\leq n (\epsd + \alpha) + \left(\frac{\lambda}{\lambda - 1}\right)^{1 - 1/\lambda} \diam(\cD)^{1/\lambda} \E\left[\left(\beta_0 + \sum_{t=1}^n \beta(\sigma_t, A_t)\right)^{1-1/\lambda}\right]\,,
\end{align*}
where the second inequality holds for any sequence of positive learning rates by elementary optimisation.
The third inequality by Fatou's lemma as in \cite[theorem 3]{LS19pminfo}.
The fourth inequality by telescoping the weighted potential and the fact that the learning rates is non-increasing.
The final inequality follows from the definition of the learning rate and standard bounding.

\paragraph{Application}
To make things concrete, let us give an application to $d$-armed bandits (see \cref{tab:examples}).
The following argument is due to \cite{BS20}.
Let $F : \R^d \to \R \cup \{\infty\}$ be the logarithmic barrier, which is defined on the positive orthant by
\begin{align*}
F(p) = - \sum_{a=1}^d \log(p_a)
\end{align*}
and is associated with Bregman divergence
\begin{align*}
\KL(p, q) = -\sum_{a=1}^d \log\left(\frac{p_a}{q_a}\right) + \ip{1/q,p-q} \,.
\end{align*}
Let $\epsilon \in (0, 1/d)$ and $\cD = \sP_\epsilon$, for which $\epsd \leq d \epsilon$. 
A simple calculation shows that $\diam(\cD) \leq d \log(1/\epsilon)$. 
Let $\beta(z, a) = z_a^2 = \Phi_a(z)^2$. The results by \cite{BS20} show that whenever $(Z, A^*)$ has law $\nu \in \cV$, then with $A$ sampled independently
from $(Z, A^*)$ with law $\E[A^*] \in \sP$, 
\begin{align*}
\E[\ip{A - A^*, \ell(Z)}] \leq \sqrt{\E[\beta(Z, A)] \E[\KL(\E[A^*|\Phi_A(Z),A], \E[A^*])]}\,.
\end{align*}
Hence, by \cref{thm:adaptive-info}, the Bayesian regret for any prior can be bounded by
\begin{align*}
\BReg_n 
&\leq nd\epsilon + \sqrt{d \E\left[\sum_{t=1}^n \ell_{A_t}(Z_t)^2\right] \log(1/\epsilon)} \\
&\leq nd\epsilon + \sqrt{d \left(\BReg_n + \E\left[\sum_{t=1}^n \ell_{A^*}(Z_t)\right]\right) \log(1/\epsilon)}
\end{align*}
Solving the quadratic shows that
\begin{align*}
\BReg_n \leq nd\epsilon + d \log(1/\epsilon) + \sqrt{d\left(1 + \E\left[\sum_{t=1}^n \ell_{A^*}(Z_t)\right]\right) \log(1/\epsilon)} \,.
\end{align*}
\cref{thm:adaptive-online} shows that a suitable instantiation of FTRL achieves about the same bound, a result which is already known \citep{LS20book}.

\section{Computation}
Given $q = Q_t \in \cD \cap \dom(\nabla F)$, \cref{alg:exp-opt} needs to compute $p \in \sPp$ and $g \in \cG$ such that
\begin{align*}
\sup_{z \in \cZ, a^* \in \cD} \Lambda_{q,\eta}(z,a^*,p, g) \leq \Lambda_\eta^* + \epsilon\,.
\end{align*}
While this is a convex optimisation problem, $\cG$ is often infinite-dimensional and the supremum need not have an explicit form.
A fundamental case where things work out is finite partial monitoring games ($\cZ$ and $\cA$ are finite).
Then all relevant quantities are finite and standard convex optimisation libraries can be used to implement \cref{alg:exp-opt} efficiently.
\cref{thm:inf} combined with the bound on the information ratio by \cite{LS19pminfo} shows that for finite non-degenerate locally observable partial monitoring games, \cref{alg:exp-opt}
enjoys a regret bounded by
\begin{align*}
\Reg_n \leq n\epsilon + 6 |\Sigma| |\cA|^{3/2} \sqrt{n \log |\cA|}\,,
\end{align*}
where $\epsilon$ is precision, which can be arbitrarily close to zero.
The same argument shows that for globally observable and (possibly degenerate) locally observable games, the algorithm also achieves the best known rates.

\section{Finite-armed bandits}
Let us now revisit the finite-armed adversarial bandit problem, which is modelled
as a linear partial monitoring game by $\cA = \{e_1,\ldots,e_d\}$, $\cZ = [0,1]^d$, $\Sigma = [0,1]$ and $\ell(z) = z$ and $\Phi_a(z) = z_a$.
\cite{AB09} used mirror descent with the standard importance-weighted estimators to design an algorithm with $\Reg_n \leq \sqrt{8dn}$, which matches the lower bound
up to constant factors \citep{ACFS95}. 
\cite{ZL19} showed that by modifying the loss estimates, mirror descent with the same potential achieves $\Reg_n \leq \sqrt{2dn} + 48d$. 
The potential function used by \cite{AB09} has the positive orthant as its domain and is defined there by 
\begin{align*}
F(q) = -2 \sum_{i=1}^d \sqrt{q_i}\,,
\end{align*}
which for $\cD = \conv(\cA)$ has $\diam(\cD) \leq 2 \sqrt{d}$.
\cite{LS19pminfo} used entropy inequalities to show that with this potential, the bandit problem has an information ratio 
of $\alpha = 0$, $\beta = \sqrt{d}$ and $\lambda = 2$. Combining this
with \cref{thm:exp-opt,thm:inf} imply that \cref{alg:exp-opt} has $\Reg_n \leq n\epsilon + \sqrt{2dn}$ for arbitrarily small $\epsilon$. 
Regrettably, however, the fact that $\Sigma$ is infinite means that the optimisation 
problem in \cref{alg:exp-opt} is infinite-dimensional.
Nevertheless, armed with the knowledge that certain loss estimation functions exist, the challenge of finding them is less daunting.
We made two guesses that made the search for a mirror descent implementation with the same bound more tractable. First, that the estimation function could be unbiased. And second,
that mirror descent with $P_t = Q_t$ would suffice. The latter guess is partially supported by \cref{thm:ts}, though here we take $\cD = \conv(\cA)$, so the conditions of the theorem
are not satisfied.
After an extended Mathematica session, an estimation function that does the job is given by
\begin{align*}
g(a, \sigma)_b =  \sind_{a=b} \left(\sigma - 1/2 + \frac{\eta}{8}\left(1 + \frac{1}{q_b + \sqrt{q_b}}\right)\right) - \frac{q_a \eta}{8(q_b + \sqrt{q_b})} \,,
\end{align*}
which is unbiased. 
Hence, mirror descent with $P_t = Q_t$ and the above estimation function has a bound on the regret of
\begin{align}
\Reg_n 
&\leq \frac{\diam(\cD)}{\eta} + \frac{1}{\eta} \E\left[\sum_{t=1}^n \Psi_{Q_t}\left(\frac{\eta g(A_t, \Phi_{A_t}(z_t))}{Q_{tA_t}}\right)\right] \nonumber \\
&\leq \frac{\diam(\cD)}{\eta} + \frac{n}{\eta} \sup_{\substack{q \in \relint(\cD) \\ z \in \cZ}} \sum_{a=1}^d q_a \Psi_{q}\left(\frac{\eta g(a, \Phi_a(z))}{q_a}\right)\,. 
\label{eq:bandit-stab} \\
&\leq \frac{\diam(\cD)}{\eta} + \frac{n \eta \sqrt{d}}{4} \nonumber \\
&= \sqrt{2nd}\,, \nonumber
\end{align}
where the final inequality follows by bounding $\diam(\cD) \leq 2 \sqrt{d}$ and choosing $\eta = \sqrt{8/n}$ and the second inequality follows from the following lemma.
Note that when $n \leq 4$, then $\Reg_n \leq \sqrt{2dn}$ is immediate. Hence we may assume that $\eta \leq \sqrt{2}$.

\begin{lemma}\label{lem:bandit-stab}
Suppose that $\eta \leq \sqrt{2}$. 
Then stability term in the right-hand side of \cref{eq:bandit-stab} is bounded by 
\begin{align*}
\frac{1}{\eta} \sup_{\substack{q \in \relint(\cD) \\ z \in \cZ}} \sum_{a=1}^d q_a \Psi_{q}\left(\frac{\eta g(a, \Phi_a(z))}{q_a}\right) \leq \frac{\eta \sqrt{d}}{4}\,.
\end{align*}
\end{lemma}

\begin{proof}
Let $z \in \cZ$ and $q \in \relint(\cD)$ be arbitrary. Then,
\begin{align*}
&\frac{1}{\eta} \sum_{a=1}^d p(a) \Psi_q\left(\frac{\eta g(a, \Phi_a(z))}{q_a}\right)
=\eta \sum_{a=1}^d q_a \sum_{b=1}^d \frac{q_b \left(\frac{g(a, \Phi_a(z))_b}{q_a}\right)^2}{\sqrt{\frac{1}{q_b}} + \frac{\eta g(a, \Phi_a(z))_b}{q_a}} \\
&\qquad=\eta \sum_{b=1}^d \sqrt{q_b} \underbracket{\left(\sum_{a=1}^d \frac{q_a \sqrt{q_b} \left(\frac{g(a, \Phi_a(z))_b}{q_a}\right)^2}{\sqrt{\frac{1}{q_b}} + \frac{\eta g(a, \Phi_a(z))_b}{q_a}}\right)}_{\textrm{(A)}_b} 
\leq \frac{\eta}{4} \sum_{b=1}^d \sqrt{q_b} 
\leq \frac{\eta \sqrt{d}}{4}\,,
\end{align*}
where the first inequality follows from the messy calculation below and the second inequality follows from Cauchy--Schwarz.
For the messy calculation:
\begin{align*}
\textrm{(A)}_b
&= \sum_{a=1}^d \frac{q_a \sqrt{q_b} \left(\frac{g(a, \Phi_a(z))_b}{q_a}\right)^2}{\sqrt{\frac{1}{q_b}} + \frac{\eta g(a, \Phi_a(z))_b}{q_a}} \\
&= \frac{1}{8}\left(\frac{\left(\eta+4 (2 z_b-1) \sqrt{q_b}\right)^2}{\eta^2+4 \eta (2 z_b-1) \sqrt{q_b}+8 q_b}+\frac{\eta^2 \left(1- \sqrt{q_b}\right)}{8 \left(\sqrt{q_b}+1\right) - \eta^2}\right) \\ 
&= \frac{1}{8}\left(2 - \frac{\eta^2}{\eta^2+4 \eta (2 z_b-1) \sqrt{q_b}+8 q_b}+\frac{\eta^2 \left(1- \sqrt{q_b}\right)}{8 \left(\sqrt{q_b}+1\right) - \eta^2}\right) \\
&\leq \frac{1}{8}\left(2 - \frac{\eta^2}{\eta^2+4 \eta \sqrt{q_b}+8 q_b}+\frac{\eta^2 \left(1- \sqrt{q_b}\right)}{8 \left(\sqrt{q_b}+1\right) - \eta^2}\right) \\
&\leq \frac{1}{4}\,,
\end{align*}
where the final inequality follows since $\eta \leq \sqrt{2}$.
\end{proof}

\section{Discussion}\label{sec:discussion}

\paragraph{Convex bandits}
Although we do not yet have an efficient approximation of \cref{alg:exp-opt} for convex bandits, the analysis here does provide some insights to that problem.
Notably, our results combined with the bound on the information ratio by \cite{Lat20-cvx} show there exist loss estimation functions and exploratory distributions 
such that \cref{alg:md} has regret at most $\Reg_n \leq O(d^{2.5} \sqrt{n} \log(n))$.
This hints towards a simpler argument than what is given by \cite{BLE17}, with no need for zooming or any sophisticated reset argument.

\paragraph{Infinite action spaces}
In principle, infinite actions spaces can be handled using the same arguments. But delicate measure-theoretic issues arise in the application of Sion's theorem and
some technical assumptions may be necessary. We leave this as a fun challenge for someone with an inclination to technical measure-theoretic details.

\bibliographystyle{plainnat}
\bibliography{all}

\appendix

\section{Technical inequalities}

Here we collect some technical results.

\begin{lemma}\label{lem:tech1}
Let $\lambda > 1$ and $(\beta_t)_{t=0}^n$ be a sequence of positive reals with $\beta_0 \geq \beta_t$ for all $1 \leq t \leq n$. Then,
\begin{align*}
\sum_{t=1}^n \beta_t \left(\sum_{s=0}^{t-1} \beta_s\right)^{1/\lambda-1} \leq \lambda \left(\sum_{t=1}^n \beta_t\right)^{1/\lambda}\,.
\end{align*}
\end{lemma}

\begin{proof}
Let $B(t) = \int_0^t \beta_{\ceil{s}} \d{s}$. Then,
\begin{align*}
\sum_{t=1}^n \beta_t \left(\sum_{s=0}^{t-1} \beta_s\right)^{1/\lambda-1}
&\leq \int_0^n B'(t) B(t)^{1/\lambda - 1} \d{t} \\
&= \lambda B(n)^{1/\lambda} \\
&= \lambda\left(\sum_{t=1}^n \beta_t\right)^{1/\lambda}\,.
\qedhere
\end{align*}
\end{proof}

\begin{lemma}\label{lem:ids}
Let $\Delta \in \R^{|\cA|}$ and $\cI \in [0,\infty)^{|\cA|}$ with $\cI \neq \zeros$ and for $p \in \sP$ let
\begin{align*}
R_\lambda(p) = \frac{\max(0, \ip{p, \Delta})^\lambda}{\ip{p, \cI}}\,.
\end{align*}
Then, for any $\lambda \geq 2$,
\begin{enumerate}
\item[(a)] $p \mapsto R_\lambda(p)$ is convex.
\item[(b)] If $p$ minimises $p \mapsto R_2(p)$, then $R_\lambda(p) \leq 2^{\lambda - 2} \min_{q \in \sP} R_\lambda(q)$.
\end{enumerate}
\end{lemma}

\begin{proof}
Part (a) follows by differentiating. For part (b), let $p$ be the minimiser of $R_2$ and $q$ the minimiser of $R_\lambda$. 
The result is immediate if $\ip{p, \Delta} \leq 0$, so assume for the remainder that $\ip{p, \Delta} > 0$.
By the first-order optimality conditions
\begin{align*}
0 \leq \ip{\nabla R_2(p), q - p} = \frac{2\ip{q - p, \Delta} \ip{p, \Delta}}{\ip{p, \cI}} - \frac{\ip{q - p, \cI} \ip{p, \Delta}^2}{\ip{p, \cI}^2}\,.
\end{align*}
Rearranging shows that
\begin{align}
\ip{p, \Delta}\left(1 + \frac{\ip{q, \cI}}{\ip{p, \cI}}\right) \leq 2 \ip{q, \Delta}\,.
\label{eq:first-order}
\end{align}
Since the information gain is non-negative, it follows that $\ip{p, \Delta} \leq 2 \ip{q, \Delta}$.
Therefore,
\begin{align*}
R_\lambda(p) 
&= \frac{\ip{p, \Delta}^\lambda}{\ip{p, \cI}} 
\leq \frac{2^{\lambda - 2} \ip{p, \Delta}^2 \ip{q, \Delta}^{\lambda - 2} }{\ip{p, \cI}}
\leq \frac{2^{\lambda - 2} \ip{q, \Delta}^2}{\ip{q, \cI}}
= 2^{\lambda - 2} \min_{q \in \sP} R_\lambda(q)\,,
\end{align*}
where the first inequality follows form \cref{eq:first-order} and the second since $p$ minimises $R_2$.
\end{proof}

The next simple lemma is used to show that the exploratory distribution can be chosen to assign non-zero probability to all actions with
arbitrarily small loss.

\begin{lemma}\label{lem:supp}
Suppose a partial monitoring game has an information ratio of $(\alpha, \beta, \lambda)$ with $\lambda \geq 1$. Then for any $\nu \in \cV$ and $\epsilon \in (0,1)$, 
there exists a $q \in \sP_\epsilon$ such that when $(Z, A^*, A)$ is sampled from the product measure $\nu \otimes q$, then
\begin{align*}
\E[\ip{A - A^*, \ell(Z)}]
\leq |\cA| \epsilon + \alpha + \beta^{1-1/\lambda} \E[\KL(\E[A^*|\Phi_A(Z), A], \E[A^*])]^{1/\lambda} \,.
\end{align*}
\end{lemma}

\begin{proof}
Let $p \in \sP$ be the distribution guaranteed by the definition of the information ratio and $q = (1 - \epsilon) p + \epsilon \ones$. Then $q \in \sP_\epsilon$, and
\begin{align*}
\sum_{a \in \cA} q(a) &\E[\ip{a - A^*, \ell(Z)}]
= (1 - \epsilon) \sum_{a \in \cA} p(a) \E[\ip{a - A^*}, \ell(Z)] \\
&\qquad\qquad\qquad\qquad\qquad\qquad + \epsilon \sum_{a \in \cA} \E[\ip{a - A^*, \ell(Z)}] \\
&\leq |\cA| \epsilon + (1 - \epsilon)\left[ \alpha + \beta^{1-1/\lambda} \left(\sum_{a=1}^k p(a) \E[\KL(\E[A^* | \Phi_a(Z)], \E[A^*])]\right)^{1/\lambda}\right] \\
&\leq |\cA| \epsilon + \alpha + \beta^{1-1/\lambda} \left( \sum_{a=1}^k q(a) \E[\KL(\E[A^*| \Phi_a(Z)], \E[A^*])]\right)^{1/\lambda} \,,
\end{align*}
where in the first inequality we used the assumption that $\ip{a, \ell(z)} \in [0,1]$ for all $a \in \cA$ and $z \in \cZ$.
The second follows by the non-negativity of the Bregman divergence and the fact that $(1 - \epsilon) \leq (1 - \epsilon)^{1/\lambda}$ since $\lambda \geq 1$ and $\epsilon \in (0,1)$.
\end{proof}

\section{Proof of \cref{thm:ts}}\label{sec:thm:ts}

Let us start with a simple lemma that, like the theorem, assumes that $F$ is the $s$-Tsallis entropy for $s \in [0,1]$.

\begin{lemma}\label{lem:b-diff}
Suppose that $\epsilon \in [-1,1]^d$ and $q \in \sP$ and $r \in [0,1]^d$, then
\begin{align*}
\ip{q - r, \epsilon} - \KL(r, q) \leq \frac{e}{2} \ip{q, \epsilon^2}\,.
\end{align*}
\end{lemma}

\begin{proof} 
It suffices to prove the result when $d = 1$. 
Let $f_s(p) = -(p^s - sp-(1-s))/(s(s-1))$, which has $f''_s(p) = p^{s-2}$.
A tedious calculation shows that the value of $r$ maximising the left-hand side satisfies $r \leq e q$. 
By Taylor's theorem and the fact that $p \mapsto f''_s(p)$ is decreasing,
\begin{align*}
\epsilon(q - r) - \KL(r, q) 
&\leq \frac{\epsilon^2}{2 f''(\max(q, r))}
= \frac{\epsilon^2}{2} (\max(q, r))^{2 - s} 
\leq \frac{e q \epsilon^2}{2}\,.
\qedhere
\end{align*}
\end{proof}

\begin{proof}[Proof of \cref{thm:ts}]
Let $\epsilon > 0$ be sufficiently small and $(Z, A^*)$ have law $\nu \in \cV_\epsilon$ and $r = \E[A^*]$.
It suffices to show that when $A$ has law $q$, then
\begin{align*}
\inf_{g \in \cG_\epsilon} \E[\Lambda(Z, A^*, p, g)] 
&= \E\left[\ip{A - A^*, \ell(Z)} - \frac{1}{\eta} \KL(\E[A^*|\Phi_A(Z), A], r) - \frac{1}{\eta} \KL(r, q)\right] \\
&\leq (1 + O(\eta^{1/2})) \frac{\eta \beta}{4}\,.
\end{align*}
Let $\cI_a = \E[\KL(\E[A^*|\Phi_a(Z)], \E[A^*])]$ and $\Delta_a = |\E[\ell_a(Z)] - \E[\ip{A^*, \ell(Z)}]|$.
Suppose first that $\ip{q, \Delta} \leq \ip{q, \cI} / \eta$. Then, by the positivity of the Bregman divergence,
\begin{align*}
\E[\Lambda(Z, A^*, q, g)]
\leq \ip{q, \Delta} - \frac{\ip{q, \cI}}{\eta} - \frac{1}{\eta} \KL(r, q) \leq 0\,.
\end{align*}
On the other hand, if $\ip{q, \Delta} > \ip{q, \cI} / \eta$, then 
\begin{align*}
&\inf_{g \in \cG_\epsilon} \E[\Lambda(Z, A^*, q, g)] 
\leq \ip{q, \Delta} - \frac{\ip{q,\cI}}{\eta} - \frac{1}{\eta} \KL(r, q) \\
&= \ip{r, \Delta} - \frac{\ip{q, \cI}}{\eta} + \ip{q - r, \Delta} - \frac{1}{\eta} \KL(r, q) \\
&\leq \sqrt{\beta \ip{r, \cI}} - \frac{\ip{q, \cI}}{\eta} + \ip{q - r, \Delta} - \frac{1}{\eta} \KL(r, q) \\
&\leq \sqrt{\beta \ip{q, \cI}} - \frac{\ip{q, \cI}}{\eta} + \frac{|\ip{r - q, \cI}| \sqrt{\beta}}{2 \sqrt{\ip{q, \cI}}} + \ip{q - r, \Delta} - \frac{1}{\eta} \KL(r, q) \\
&\stackrel{(\star)}\leq \sqrt{\beta \ip{q, \cI}} - \frac{\ip{q, \cI}}{\eta} + \eta \bip{q, \left(\frac{\cI \sqrt{\beta}}{\sqrt{q, \cI}} + 2\Delta\right)^2} - \frac{1}{2\eta} \KL(r, q) \\
&\leq \sqrt{\beta \ip{q, \cI}} - \frac{\ip{q, \cI}}{\eta} + \frac{2\eta \beta}{\ip{q, \cI}} \bip{q, \cI^2} + 8\eta \ip{q, \Delta^2} - \frac{1}{2\eta} \KL(r, q) \\
&\leq \sqrt{\beta \ip{q, \cI}} - \frac{\ip{q, \cI}}{\eta} + \frac{2\beta \ip{q, \cI}}{\eta^{1/3}} + 8\eta \ip{q, \Delta} - \frac{1}{2\eta} \KL(r, q) \\
&\leq \sqrt{\beta \ip{q, \cI}} - \left(1 - 2\beta \eta^{2/3} - 8\eta \right) \frac{\ip{q, \cI}}{\eta} + 8\eta \left(\ip{q, \Delta} - \frac{\ip{q, \cI}}{\eta} - \frac{1}{\eta} \KL(r, q)\right)\,,
\end{align*}
where the first inequality follows from \cref{eq:dual}, 
the second by assumption and the third since $(x + \delta)^{1/2} \leq x^{1/2} + \frac{1}{2} |\delta| x^{-1/2}$.
The fifth inequality follows from the fact that $(x+y)^2 \leq 2x^2 + 2y^2$ and the sixth since $\ip{q, \cI^2} \leq \ip{q, \cI}^2 / \eta^{4/3}$ and $\Delta \leq \ones$.
The last inequality follows from naive simplification and re-arranging and by taking $\eta$ suitably small.
The inequality marked with a $(\star)$ follows from \cref{lem:b-diff}, which is justified because
\begin{align*}
\frac{\eta \cI_a \sqrt{\beta}}{\sqrt{\ip{q, \cI}}} + 2\eta \Delta_a
&\leq \frac{\eta \sqrt{\beta \ip{q, \cI}}}{q_a} + 2\eta \Delta_a 
\leq \frac{\eta \sqrt{\beta \eta}}{q_a} + 2\eta \Delta_a 
\leq 1\,,
\end{align*}
which holds for all sufficiently small $\eta$ since $q_a \geq \eta^{4/3}$.
Rearranging shows that
\begin{align*}
\E[\Lambda(Z, A^*, q, g)]
&\leq \frac{1}{1 - 8\eta} \left(\sqrt{\beta \ip{q, \cI}} - \left(1 - 2\beta \eta^{2/3} - 8\eta\right) \frac{\ip{q, \cI}}{\eta}\right) \\
&= (1 + O(\eta^{2/3})) \frac{\eta \beta}{4}\,.
\end{align*}
\end{proof}

\end{document}